\documentclass[11pt]{amsart}

\usepackage[english]{babel}

\usepackage{amsmath,amsfonts,amsthm}
\def\ds{\displaystyle}
\def\R{{\mathbb R}}
\def\Z{{\mathbb Z}}

\def\N{{\mathbb N}}
\def\ali{\hfill\break}

\newtheorem{defi}{Definition}
\newtheorem{theo}{Theorem}

\newtheorem{cor}{Corollary}

\author[N. Enriquez]{Nathana\"el ENRIQUEZ}
\address{Laboratoire de Probabilit\'es et mod\`eles al\'eatoires,
4 place Jussieu,
75252 Paris cedex 05}
\email{enriquez@ccr.jussieu.fr}

\author[C. Sabot]{ Christophe SABOT}
\address{Laboratoire de Probabilit\'es et mod\`eles al\'eatoires,
Universit\'e Paris 6, 4
place Jussieu,
75252 Paris cedex 05}
\email{sabot@ccr.jussieu.fr}

\title[Reinforced random walks]{A note on
edge oriented reinforced
random walks and RWRE.}

\begin{document}

\begin{abstract}
  This work introduces the notion of edge oriented reinforced random
walk which proposes in a
general framework an alternative understanding of the annealed
law of random walks in random environment.
\end{abstract}
\maketitle
\section{Introduction}
The aim of this text is to present in a general framework
a precise and explicit
correspondance between a class of edge oriented reinforced random 
walks and random walks in random environment (RWRE).
In other words the study of an
edge oriented reinforced random walks, which is strongly non-Markovian,
is equivalent by introducing an external randomness, to a 
Markovian problem.
This relation already appeared, as a tool, in a work of Pemantle
(cf \cite{Pemantle}),
for the study of reinforced random walks
on binary trees, for very special laws of reinforcement.

The essence of the result is that 
the law of an edge oriented reinforced random walk
coincides with the annealed laws
of a RWRE
if the law of reinforcement satisfies the following condition:
at one point the probability of a sequence of successive moves
depends only on the number of each type of move. 
This condition, we call admissiblity, is expressed as the closedness of a 
certain discrete form.
The proof of our result is based on a theorem on the moment problem
for random variables in $[0,1]^d$.
As a corrolary, the uniqueness of this moment problem implies
that the annealed law determines the quenched law for RWRE.

\section{Definitions and statement of the result}
\begin{defi}
We call law of reinforcement 
with $d$ neighbours a function $V$:
$$V:\Z_+^d\mapsto T_d:=\{(x_1,...,x_d)\in ]0,1]^d, \,
s.t.,\,\sum_{i=1}^dx_i=1\}$$
$$\vec{p}=(p_1,...,p_d)\to(V_1(\vec{p}),...,V_d(\vec{p}))$$
\end{defi}

We now define for a law of reinforcement the notion of admissibility, which
will play a key role in the following. Let $(e_i)_{1\leq i\leq d}$ denote the
canonical basis of $\Z_+^d$.
\begin{defi}
Let us consider the graph $\Z_+^d$ whose edges $(\vec{p},\vec{p}+e_i)$ are
oriented from $\vec{p}$ to $\vec{p}+e_i$. Let $V$ be a law of reinforcement
on $\Z_+^d$ and $w$ be the 1-form defined by
$w((\vec{p},\vec{p}+e_i))=\ln V_i(\vec{p})$.

We shall say that $V$ is admissible when $w$ is a closed form.

\end{defi}

We consider now a countable graph $G$ having at any point a finite number
of neighbours. For any vertex $x$ of $G$ we denote by $d(x)$ the cardinal
of the neighbours of $x$ and $\{e(x,i),i=1,...,d(x)\}$ the neighbours of
$x$. At any vertex we suppose
given a law of reinforcement $V(x)$ with $d(x)$
neighbours: $V(x) : \vec{p}\to(V_i(x,\vec{p}))_{1\leq i\leq d(x)}$.

\begin{defi}
We call reinforced random walk with law of reinforcement  $V(x)$, the random walk
defined by the family of laws on the trajectories starting at $x_0$,
$(\hat{P}_{x_0})_{x_0\in G}$ given by
$$\hat{P}_{x_0}(X_{n+1}=e(x,i)|\sigma(X_n=x)\wedge\sigma (X_k,k\leq
n-1))=V_i(x,\vec{N}(x))$$ where 
$\vec{N}(x)=(N_i(x))_{1\leq i\leq d(x)}$ and
$N_i(x)=\sum_{l=0}^{n-1}1_{\{X_l=x,X_{l+1}=e(x,i)\}}$.

\end{defi}

\begin{defi}
A reinforced random walk is called admissible when $V(x)$ is admissible for all
vertex $x$ of $G$.
\end{defi}

Let us now introduce  random walks in random environment on the graph $G$.

We define an environment as an element
$\omega=(\omega(x))_{x\in G}$ where at any
vertex $x$, $\omega(x)$ is in $T_{d(x)}$.
At any vertex $x$ of $G$, we consider a probability measure $\mu_x$ on
$T_{d(x)}$ and we set $\mu:=\ds\mathop{\otimes}_{x\in G}\mu_x$,
so that $\mu$ is a probability measure on the environments such that
$(\omega(x))_{x\in G}$ are independent
random variables of law $\mu_x$.

We denote by $P_{x_0,\omega}$ the law of the Markov chain in the
environment $\omega$ starting at $x_0$ defined by:
$$\forall x_0\in G,\;\; \forall k\in\N,\;\;
P_{x_0,\omega}(X_{k+1}=e(x,i)|X_k=x)=\omega(x,i).$$
Finally we denote by $P_x$ the annealed measure i.e.
$P_x=\mu\otimes P_{x,\omega}$.
We are now able to state our main result:

\begin{theo}
For any countable graph $G$

i) For all law  of environment $\mu=\ds\mathop{\otimes}_{x\in G}\mu_x$
the law of the reinforced random walk $(\hat{P}_{x_0})_{x_0\in G}$
associated with the law of reinforcement $V$
given by:
\begin{eqnarray}
\label{1}
V_i(x,p_1,...,p_{d(x)})=
{E_{\mu_x}[\omega(x,i)\prod_{j=1}^{d(x)}\omega(x,j)^{p_j}]\over
E_{\mu_x}[\prod_{j=1}^{d(x)}\omega(x,j)^{p_j}]}
\end{eqnarray}
(where $E_{\mu_x}$ denotes the expectation under the law $\mu_x$)
coincides with the annealed law of the RWRE $(P_{x_0})_{x_0\in G}$.
Moreover the law $(V(x))$ is admissible.

ii) Conversely, if $V$ is an admissible law of reinforcement on $G$, then
there exists a unique law of environment $\mu=\ds\mathop{\otimes}_{x\in G}\mu_x$
for which equality (\ref{1}) is satisfied, thus for which the law of the reinforcement 
random walk 
$(\hat{P}_{x_0})_{x_0\in G}$ 
coincides with the annealed law $(P_{x_0})_{x_0\in G}$.
\end{theo}
\begin{cor}
The annealed law determines the quenched law for RWRE, i.e. 
there can not exist two different laws of environment $\mu$ and
$\mu'$ with the same annealed law $(P_{x_0})_{x_0\in G}$. 
\end{cor}
Example 1: At any point $x$, choose a vector
$(\alpha(x,1),...,\alpha(x,d(x))$ in $(\R_+^*)^{d(x)}$. The law of reinforcement 
$V_i(x,p_1,..., p_{d(x)})={\alpha(x,i)+p_i\over
\sum_{j=0}^{d(x)}\alpha(x,j)+p_j}$ 
is an admissible law of reinforcement
associated with the environment $(\mu_x)_{x\in G}$ where $\mu_x$ is a
Dirichlet law with parameters $(\alpha(x,1),...,\alpha(x,d(x)))$, i.e. $\mu_x$ is
the law on $T_{d(x)}$ with density
$$\mu_x(t_1,...,t_{d(x)})=
{\Gamma(\alpha(x,1)+...+\alpha(x,d(x))
\over
\prod_{i=1}^{d(x)}\Gamma(\alpha(x,i))}\prod_{i=1}^{d(x)}t_i^{\alpha(x,i)-1}.$$
This is the meaning of the classical Polya's urn scheme
(cf \cite{Feller}, section VI.12), if we put at
all sites an independent urn which will define the move at that site.
\ali
Example 2:
We can generalize the previous example as follows:
at all point $x$ in $G$, choose not only a vector $(\alpha(x,1),...,\alpha(x,d(x))\in (\R_+^*)^{d(x)}$, but
also an integer $n(x)$ and a homogeneous polynomial 
$P(x, t_1,\ldots , t_{d(x)})$ of degree $n(x)$ 
of the form:
$$P(x, t_1,\ldots , t_{d(x)})=
\sum_{{k_1,\ldots ,k_{d(x)},
\atop
k_1+\cdots +k_{d(x)}=n(x)}}
a_{k_1,\ldots ,k_{d(x)}}(x) t_1^{k_1} \cdots t_{d(x)}^{k_{d(x)}},
$$
where the $a_{k_1,\ldots ,k_{d(x)}}(x)$'s are non-negative reals 
(not all null).
Then we define $\mu_x$ as the law on $T_{d(x)}$ with density
$$ \mu_x(t_1,...,t_{d(x)})=
{\left(
\prod_{i=1}^{d(x)}t_i^{\alpha(x,i)-1} 
\right)
P(x,t_1,\ldots , t_{d(x)})
\over
\int_{T_{d(x)}}
\left( 
\prod_{i=1}^{d(x)}t_i^{\alpha(x,i)-1}
\right) 
P(x,t_1,\ldots , t_{d(x)})
}.
$$
On the other side, consider the polynomials
$$Q(x, y_1,\ldots ,y_{d(x)})=
\sum_{{k_1,\ldots ,k_{d(x)},
\atop
k_1+\cdots +k_{d(x)}=n(x)}}
a_{k_1,\ldots ,k_{d(x)}}(x)
\prod_{i=1}^{d(x)} (y_i,k_i),
$$
where we  write $(y,k)$ for the product $y\cdots (y+k-1)$.
Then the law $(\mu_x)$ is associated with the law of reinforcement
$V(x)$, where $V_i(x,p_1,\ldots ,p_{d(x)})$ is given by (to simplify, we forget
the $x$ dependance in the next formula, and simply write $\alpha_i$ for $\alpha(x,i)$ and $n$ for $n(x)$)
$$
\left(
{\alpha_i+p_i\over (\sum_{j=1}^d \alpha_j +p_j) +n}
\right)
{ Q(\alpha_1 +p_1,\ldots , \alpha_i+p_i+1, \ldots ,\alpha_d+p_d)
\over
Q(\alpha_1 +p_1,\ldots ,\alpha_d+p_d)}.
$$

  Proof:

i) For any vertices $x_0,x$ of $G$,
$\forall 1\leq i\leq d(x), \forall n\in\N,$
\begin{eqnarray*}
&&P_{x_0}(X_{n+1}=e(x,i)|(X_n=x)\wedge\sigma(X_k,k\leq n-1))
\\
%$$={P_{x_0}(X_{n+1}=e(x,i),(X_n=x)\wedge\sigma(X_k,k\leq n-1)\over
%P_{x_0}((X_n=x)\wedge\sigma(X_k,k\leq n-1)}$$
&=&{E[\omega(x,i)\prod_{y\in
G}\prod_{j=1}^{d(y)}\omega(y,j)^{N_j(y)}]\over
E[\prod_{y\in G} \prod_{j=1}^{d(y)}\omega(y,j)^{N_j(y)}]},
\end{eqnarray*}
where $N_j(y)$ is as defined in definition 3. Now using the independance of
the variables $\omega(y,i)$ for different vertices $y$, the terms depending
on $\omega(y,j)$ for $y\neq x$ cancel in the previous ratio  and we
get i).

ii) The only point is to prove that for any admissible law of reinforcement
$V$ with $d$ neigbours 
there exists  a	probability measure
$\mu$ on $T_d$, such that a $T_d$-valued random variable
$\vec{X}:=(X_1,...,X_d)$ of law $\mu$  satisfies
$V(p_1,...,p_d)={{E[X_i\prod_{j=1}^d X_j^{p_j}]\over
E[\prod_{j=1}^d X_j^{p_j}]}}$.

We can see that the condition is expressed in terms of the moments of
$\mu$, therefore the aim is to prove that the assumption of
admissibility on $V$ implies the solvability of the moment problem for
$\mu$.

For that purpose we introduce the quantities which are intended to be the
moments of $\mu$.

Consider $s\in\{1,...,d\}^\N$ and the path $U(s,.)$ on $\Z_+^d$
defined by:

$U(s,n):=\sum_{i=0}^{n-1}e_{s(i)}$

and let $M(s,n):=\prod_{i=0}^{n-1}V_{s(i)}(U(s,i))$.

The fact that $V$ is admissible implies that $M(s,n)$ depends only on
$U(s,n)$, indeed $M(s,n)=\exp(\int_0^{U(s,n)} w)$ (with the notations
of definition 2). So let us introduce the sequence with $d$ indices:
$$v_{k_1,...,k_{d}}:=M(s,n) \hbox{ for }
U(s,n)=k_1e_1+...+k_{d}e_{d}$$
  What remains to prove is that $v$ is the sequence of moments of a
$T_{d(x)}$-valued variable. For that purpose we use the generalization of
the Hausdorff criterium, concerning the existence of a (unique) solution
to the moment problem for random variables on
$[0,1]^{d}$. For any $\vec{h}={h_1,\ldots ,h_d}$ we
define the operator $\Delta^{\vec{h}}$ on real sequences indexed
by $\Z_+^d$, i.e. $\Delta^{\vec{h}}:\R^{\Z_+^d}\rightarrow \R^{\Z_+^d}$
and defined recursively by
$$\Delta^{e_i}(u)=(u_{\vec{k}+e_i}-u_{\vec{k}})_{\vec{k}\in \Z_+^d},
$$
and 
$$\Delta^{\vec{k}+e_i}=\Delta^{e_i}\circ \Delta^{\vec{k}}.
$$
(Remark that this definition is valid since the $\Delta^{e_i}$'s commute).
We recall here the result of Hildebrandt and Schoenberg
(cf \cite{Hildebrandt}): a sequence 
$(u_{\vec{k}})\in [0,1]^{\Z_+^d}$
is the moment sequence of a probability measure on $[0,1]^d$,
i.e. $u_{\vec{k}}= \int t_1^{k_1}\cdots t_d^{k_d} d\mu(t_1,\ldots ,t_d)$
if and only if for all $\vec{h}$ and $\vec{k}$ in $\Z_+^d$, 
$(-1)^{\sum h_i} \Delta^{\vec{h}} (u)(\vec{k})$ is
positive.
\ali
Let us verify this for the sequence $v_{\vec{k}}$ introduced
previously.
Since for all $\vec{k}$, $\sum_{i=1}^d V_i(\vec{k})=1$
we have:
$$-\Delta^{e_i}(v)(\vec{k})=
\sum_{{j=1\atop j\neq i}}^d v_{k_1,\ldots, k_j+1,\ldots ,k_d}.$$
Hence, by composition
$(-1)^{h_1+...+h_{d}}\Delta^{h_1,...,h_{d}}(v_{k_1,...,k_{d}})$
will always be a linear combination with positive coefficients
of the terms of the sequence $v$. So
the condition of the criterium is satisfied.
Hence $\mu$ exists
and is unique as a solution of the moment problem
whose support is compact.
The last thing to check is that it is supported by
$T_{d}$. Let $\vec{X}$ be a random variable with law $\mu$.
The only thing to check is that $\sum_{i=1}^d X_i=1$ $\mu$-almost surely.
Since the law of $\sum_{i=1}^d X_i$ has compact support this is equivalent to
show that all the moments of $\sum_{i=1}^d X_i$ are equal to 1.
Using the fact that $\sum_{i=1}^d V_i(\vec{k})=1$ at all point
$\vec{k}$ we know that for all integer $n$
\begin{eqnarray*}
1&=&
\sum_{s\in\{1,\ldots ,d\}^n} M(s,n)\\
&=& \sum_{k_1+\cdots +k_d=n}
\#\{s\in\{1,\ldots,d\}^n,\; U(s,n)=(k_1,\ldots ,k_d)\}
v_{k_1,\ldots ,k_d}
\\
&=&
\sum_{k_1+\ldots +k_d=n} 
{(k_1+\cdots +k_d)!\over k_1!\cdots k_d!} v_{k_1,\ldots ,k_d},
\end{eqnarray*}
and this last expression is the $n$-th moment of $(X_1+\ldots +X_d)$.


\begin{thebibliography}{10}
\bibitem{Feller}
Feller, W.,
 An Introduction to Probability Theory and Its Applications.
Vol. II. John Wiley \& Sons, Inc., New York, N.Y., 1950.
\bibitem{Hildebrandt}
Hildebrandt, T. H., Schoenberg, I. J.,
{\it
On linear functional operators  and the moment problem
for a finite interval in one or several dimensions}.
Annals of Math (2), 1933, vol. 34, 317-328. 
\bibitem{Pemantle}
Pemantle, R.,
{\it
Phase transition in reinforced random walk and RWRE on trees.}
 Ann. Probab. 16 (1988), no. 3, 1229--1241.
\end{thebibliography}
\end{document}